\documentclass[smallextended,numbook,runningheads]{svjour3}     
\smartqed  
\usepackage{amsmath, amsthm, amssymb, amsfonts, enumerate}
\usepackage{graphicx}
\usepackage{float}
\usepackage{hyperref}

\usepackage{mathptmx}      
\usepackage{epsfig,color}
\usepackage{subcaption}
\usepackage[labelfont=bf]{caption}

\begin{document}
	\title{Blaschke and Separation Theorems 
		for Orthogonally Convex Sets
	}
	
		\titlerunning{ \small Blaschke and separation theorems 
			for orthogonally convex sets}        

	\author{Phan Thanh An \and Nguyen Thi Le}
		
	\institute{P. T. An\at
			Institute of Mathematical and Computational Sciences (IMACS)\\
			Faculty of Applied Science, 
			Ho Chi Minh City University of Technology (HCMUT), 268 Ly Thuong Kiet Street, District 10, Ho Chi Minh City, Vietnam\\
			Vietnam National University Ho Chi Minh City, Linh Trung Ward, Thu Duc City, Ho Chi Minh City, Vietnam \\
			\email{thanhan@hcmut.edu.vn}
			\and
		Nguyen Thi Le 
		\at
		Institute of Mathematics,  Vietnam Academy of Science and Technology (VAST), 18 Hoang Quoc Viet Street, 
		Hanoi, Vietnam
	}
	
	
	
\date{Received: date / Accepted: date}

\maketitle
	
	\begin{abstract}
In this paper, we 
 establish a Blaschke-type  theorem for path-connected and orthogonally convex sets in the plane using orthogonally convex paths.  The separation
	of these sets is established  using suitable grids.  Consequently, a closed and orthogonally convex set is represented by the intersection of staircase-halfplanes in the plane. Some topological properties of orthogonally convex sets in 
	finite-dimensional spaces  are also given.
	
	\end{abstract}
	

	\section{Introduction}\label{Introduction}%

The study of convex sets plays a significant role in geometry, analysis, and linear algebra.
 It is natural to relate convexity to orthogonally convexity (ortho-convexity, for briefly).  An ortho-convex set is also known as a rectilinear, $x$-$y$ convex, or separate convex set which is a special case of $D$-convex sets and restrict oriented convex sets~\cite{Fink2012,Matousek1998,Nicholl1983,Ottmann1984,RawlinsThesis1987}.
 Ortho-convexity has found applications in several  research  fields,  including  digital images processing~\cite{Aman2020}, VLSI circuit layout design~\cite{Biedl2011}, and geometric search~\cite{AnHai2022,SonHwangAhn2011}.
 
 The notation of ortho-convexity has been presented since the early eighties \cite{Unger1959}.
 A set in the plane is said to be ortho-convex if it intersects with  every
 horizontal or vertical line that is connected. 
  Thus every convex set is ortho-convex but not vice versa.
After that, the notion of ortho-convexity is extended to  higher dimensional spaces~\cite{Fink1998}.
Researchers have studied and given some properties of ortho-convex
sets  along with developing results similar to that of
traditional convex sets~\cite{Alegria2021,AnHuyenLe2021,Fink2012,RawlinsThesis1987}. 
 However, most of these results are obtained in the plane.
 %
  Some algorithms for finding the orthogonally convex hull of a set of finite points or of $x$-$y$ polygons  are found in~\cite{Aman2020,AnHuyenLe2021,Franek2009,Gonzalez2019,Linh2022,Nicholl1983}, where an $x$-$y$ polygon is a simple polygon with horizontal and vertical
  edges.

For  nontraditional convex sets, classical theorems for convex sets are normally considered to apply, such as Motzkin-type theorem, Radon-type theorem, Helly-type theorem, 
Blaschke-type theorem, and the separation for geodesic convex sets~\cite{AnGiangHai2010,Hai2011}, Helly-type theorem for roughly convexlike sets~\cite{An2007}.
About ortho-convex sets,
Matou{\v{s}}ek and Plech{\'a}{\v{c}} in
\cite{Matousek1998} presented Krein-Milman-type theorem for functional $D$-convex sets which belong to a class of ortho-convex sets.
Krasnosel'skii theorem for staircase path in orthogonal polygons in the plane was addressed by Breen in~\cite{Breen1994}. A natural question is “whether a Blaschke-type theorem holds for ortho-convex sets?". 

The
 separation property of ortho-convex sets in the plane has been established in many previous works in~\cite{Dulliev2017,Fink1998,Fink2012,RawlinsWood1991}. 
 For a point outside of a path-connected, closed, and ortho-convex set, there exists a so-called ortho-halfplane (see~\cite{Fink1998}) containing this point and separated from that set.
   The separation of 
   two connected components of an ortho-convex set is thus obtained (see~\cite{RawlinsWood1991}). Then the
 separation of two disjoint ortho-convex sets has been shown by Dulliev~\cite{Dulliev2017}, where the boundary of an ortho-halfplane separating two disjoint ortho-convex sets is an ortho-convex path which is a path being ortho-convex.

 In this paper, we establish a Blaschke-type theorem for path-connected and ortho-convex sets in the plane  using ortho-convex paths (Theorem~\ref{theo:Blaschke}). 
The separation of such sets is also shown (Theorem~\ref{theo:separate-2sets}) which partly answers for Conjecture 7.4 in~\cite{Fink2012} using suitable grids.
Herein, the concept of a halfplane corresponding to ortho-convexity is a staircase-halfplane whose boundary is a staircase line (see Definition~\ref{defi:staircase-line} for more details). Note that a staircase line is the specified case of an ortho-convex path. Thus the separation of two disjoint ortho-convex sets by a staircase line implies the separation of two disjoint ortho-convex sets by an ortho-convex path.
Consequently, a closed, path-connected, and ortho-convex set in the plane can be represented by the intersection of  staircase-halfplanes (Theorem~\ref{theo:represent-close-o-conv}).
Some topological properties of ortho-convex sets in $n$-dimensional space $\mathbb{R}^n$ are also presented.

	\section{Preliminary}

		We now recall some basic concepts and properties. 
		For any points $a,b$ in $\mathbb{R}^n$, we denote 
		$[a,b] := \{ (1-\lambda)a + \lambda b: 0 \leq \lambda \leq 1 \}$, 
		$]a,b] := [a,b] \setminus \{ a \}$, 
		$]a,b[ := [a,b] \setminus \{ a,b \}.$ 
		Let us denote $B(a,r) = \{ p \in \mathbb{R}^n | \, \| p-a \| <r \}$, 
		where $r>0$, $\|a\|$ is Euclidean norm, i.e., for $a=(a_i)_{i=1}^n$,
		$ \|a\| = \sqrt{\sum_{i=1}^{n}a_i^2}$.
A \emph{path} in $S \subset \mathbb{R}^n$ is a continuous mapping $\gamma$ from
an interval $I\subset  \mathbb{R}$ to $S$, $\gamma: I \to S$. If $I=[t, t']$, we say that $\gamma$ \emph{joins}
the point $\gamma(t)$ to the point $\gamma(t')$, and then the \emph{length} of $\gamma$ is the quantity
$
length(\gamma)=\sup_{\sigma} \sum_{i=0}^{k-1} \|\gamma(t_{i})-\gamma(t_{i+1})\|,
$
where the supremum is taken over the set of partitions $t=t_0<t_1<\cdots<t_k=t'$
of $[t, t']$. 
%
By abuse of notation, sometimes we also call the image $\gamma([t, t'])$  the path $\gamma: [t, t']\to S$.
For all $u,v\in \gamma, u=(u_x,u_y), v=(v_x,v_y)$, if $(u_x -v_x)(u_y -v_y) \le 0$ (or $\ge 0$), then $\gamma$ is said to be \textit{$xy$-monotone}.

Given two nonempty sets  $A$ and $B$ in $\mathbb{R}^n$, we set $$d(A,B)=\inf \{\|a-b\|, a \in A, b \in B\}$$
and $d(x, A) := d(\{x\}, A) := \inf_{a \in A}\|x-a\|$ for $x \in \mathbb{R}^n$. If $A$ and $B$ are nonempty
and compact, then $A \cap B = \emptyset$ if and only if $d(A, B) > 0.$
 \textit{Hausdorff distance} between $A$ and $B$, denoted by $d_{\mathcal{H}}(A,B)$ is defined as $$d_{\mathcal{H}}(A,B)=\max \{ \sup_{x \in A} d(x,B), \sup_{y \in B} d(y,A) \}.$$
 Then $d_\mathcal{H}$ is a metric over the space of nonempty compact subsets of $\mathbb{R}^n$.
Most of basic topological properties in metric spaces such as openness, closeness, connectedness, path-connectedness, closure, and interior operations of a subset $S$ ($\text{cl}S$ and $\text{int}S$)   can be found in~\cite{Lay2007}.

	A set in $\mathbb{R}^n$ is \textit{convex}
	if the line segment joining any two points inside it lies completely inside the set.	Then the intersection between a convex set and any hyperplane is convex.  In $\mathbb{R}^1$, an  \textit{ortho-convex} set and a convex set are the same, which is a connected set (an empty set, a point, or a single interval).
In $\mathbb{R}^2$, a set is \textit{ortho-convex} 
	if its intersection with any horizontal or vertical line is connected (see  \cite{AnHuyenLe2021,Unger1959}). Fig.~\ref{fig:examp-non-seperate}(i) gives an example of an ortho-convex set $S$. 
	A line (line segment, resp.) parallel to one of the coordinates axes is called an \textit{axis-aligned line} (\textit{axis-aligned line segment}, resp.). A rectangle having four edges that are axis-aligned line segments is called an \textit{axis-aligned rectangle}.
	Next part, we introduce the notion of staircase-halfplane in $\mathbb{R}^2$.

		\begin{definition}[\cite{Breen1994}]
		In $\mathbb{R}^2$, a \textit{staircase segment} is a simple polyline $v_0v_1 \ldots v_n$ formed by axis-aligned segments and for $i$ odd, the vectors $\overrightarrow{v_{i-1}v_i}$ have the same direction, and for $i$ even,  the vectors $\overrightarrow{v_{i-1}v_i}$ have the same direction, $1 \le i \le n$.
	\end{definition}
	
\begin{figure}
	\centering
	\begin{subfigure}[b]{0.3\textwidth}
		\centering
		\includegraphics[scale=0.2]{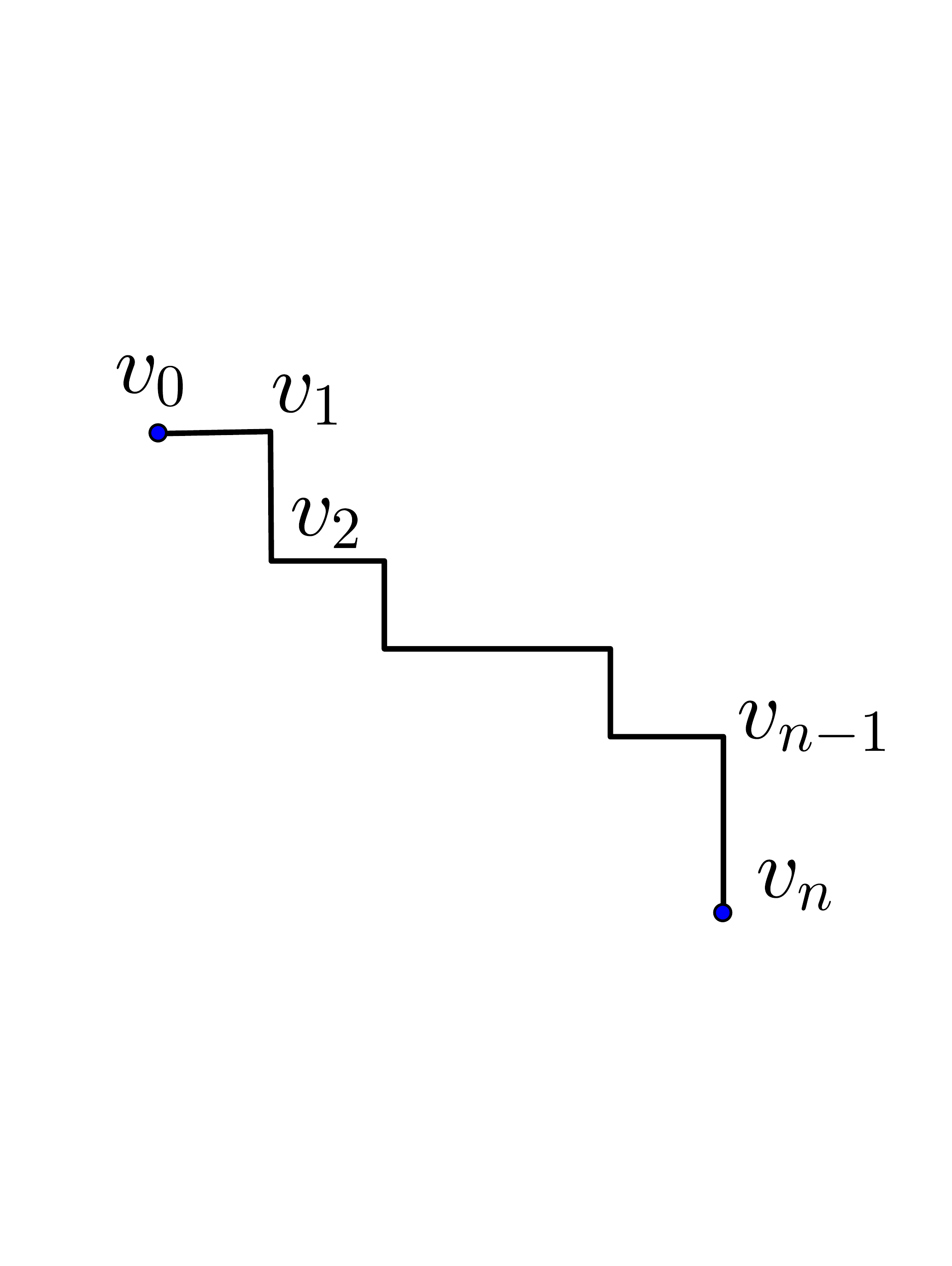}
		\caption*{(i)}
	\end{subfigure}
	\begin{subfigure}[b]{0.33\textwidth}
		\centering
		\includegraphics[scale=0.2]{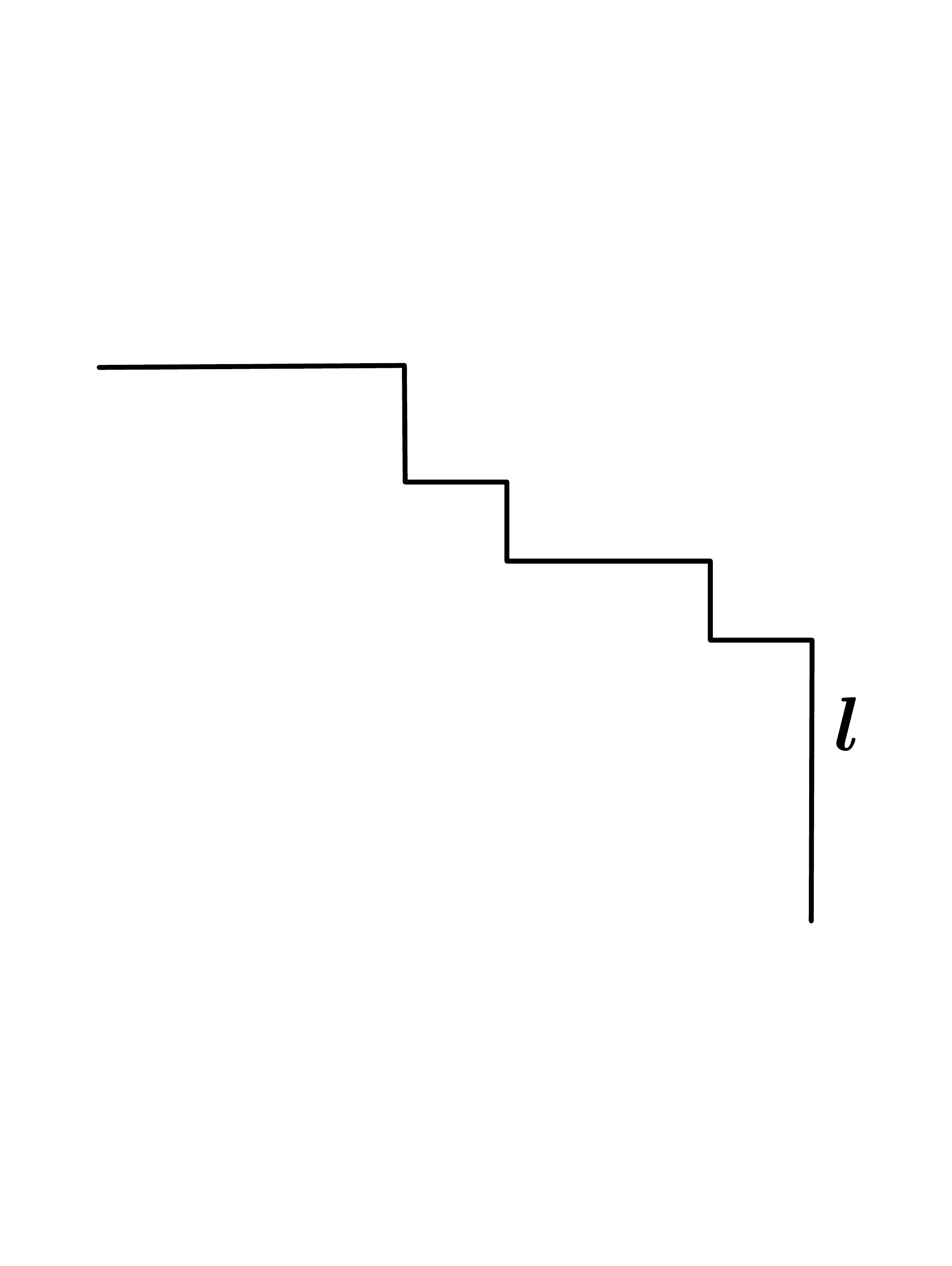}
		\caption*{(ii)}
	\end{subfigure}
	\begin{subfigure}[b]{0.33\textwidth}
		\centering
		\includegraphics[scale=0.2]{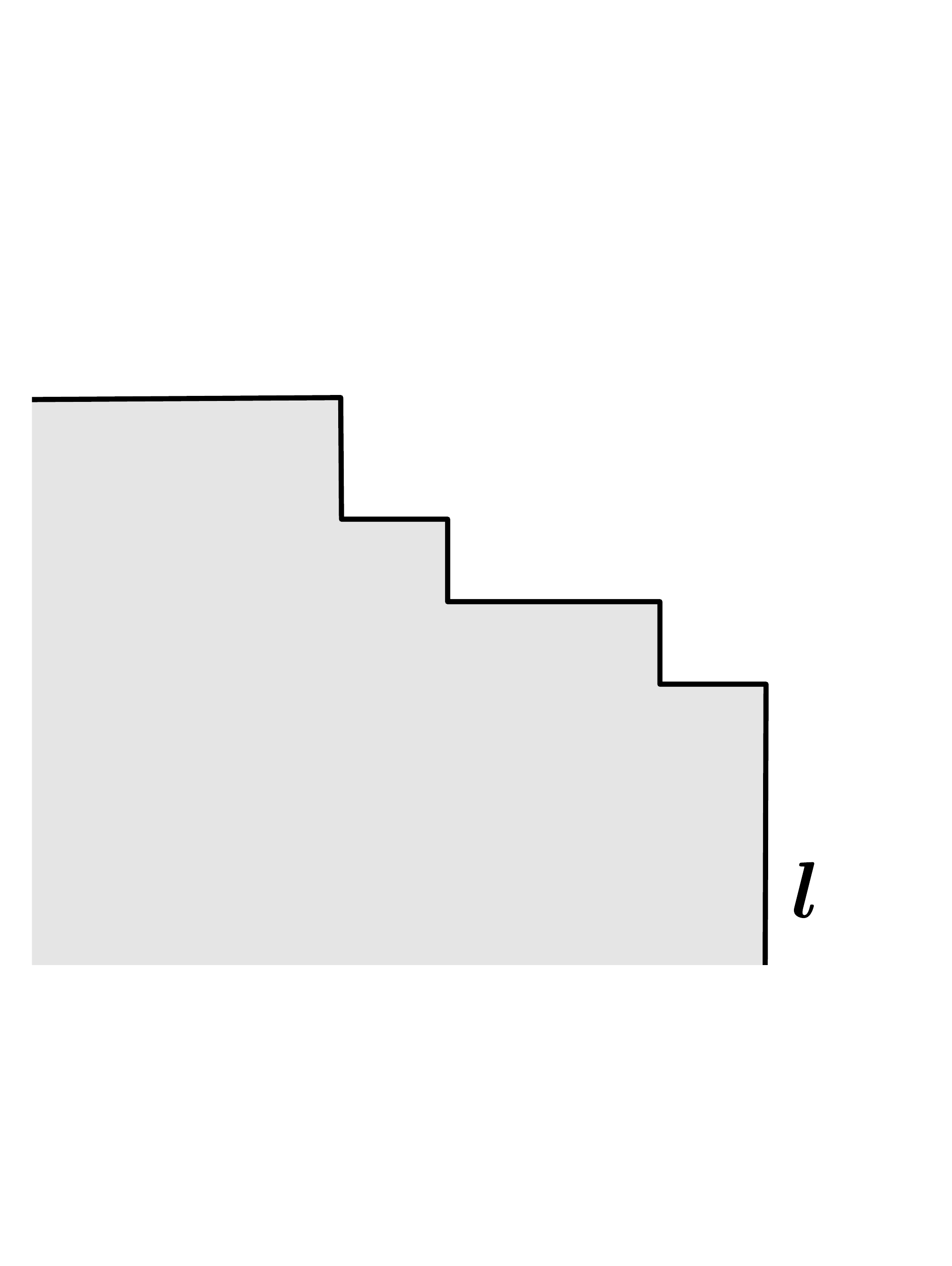}
		\caption*{(iii)}	
	\end{subfigure}
	\caption[]{Illustration of a staircase segment $v_0v_1 \ldots v_n$, the staircase line $l$ formed by the staircase segment, and a staircase-halfplane formed by $l$.}
	\label{fig:staircase_objects}
\end{figure}

	\begin{definition}
		\label{defi:staircase-line}
	In $\mathbb{R}^2$, a	\textit{staircase line} is a path obtained from a staircase segment $v_0v_1 \ldots v_n$ by replacing line segments $v_0v_1$ and $v_{n-1}v_{n}$ with the rays $v_1v_0$ and $v_{n-1}v_n$, respectively.
		
	\end{definition}
	Figs.~\ref{fig:staircase_objects}(i) and (ii) show an example of a staircase segment $v_0v_1 \ldots v_n$ and the staircase line $l$ formed by the staircase segment.
	Clearly, staircase segments and staircase lines are ortho-convex. Furthermore, axis-aligned segments (reps. axis-aligned lines) are staircase segments (reps. staircase lines).
	With respect to ortho-convex sets,
	staircase segments (reps. staircase lines)  play the same role as line segments (reps. straight lines) do for convex sets.
	 We give a new concept of a so-called staircase-halfplane relative to ortho-convexity whose boundaries made up of staircase lines (Definition~\ref{def:staircase-halfplane} and Fig.~\ref{fig:staircase_objects}(iii)). They are  analogs of half-planes.
	
	\begin{definition}
		\label{def:staircase-halfplane}
		In $\mathbb{R}^2$, a \textit{staircase-halfplane} is a closed set whose boundary is a staircase line.
	\end{definition}

\begin{figure}[htp]
		\centering
		\includegraphics[width=0.66\linewidth]{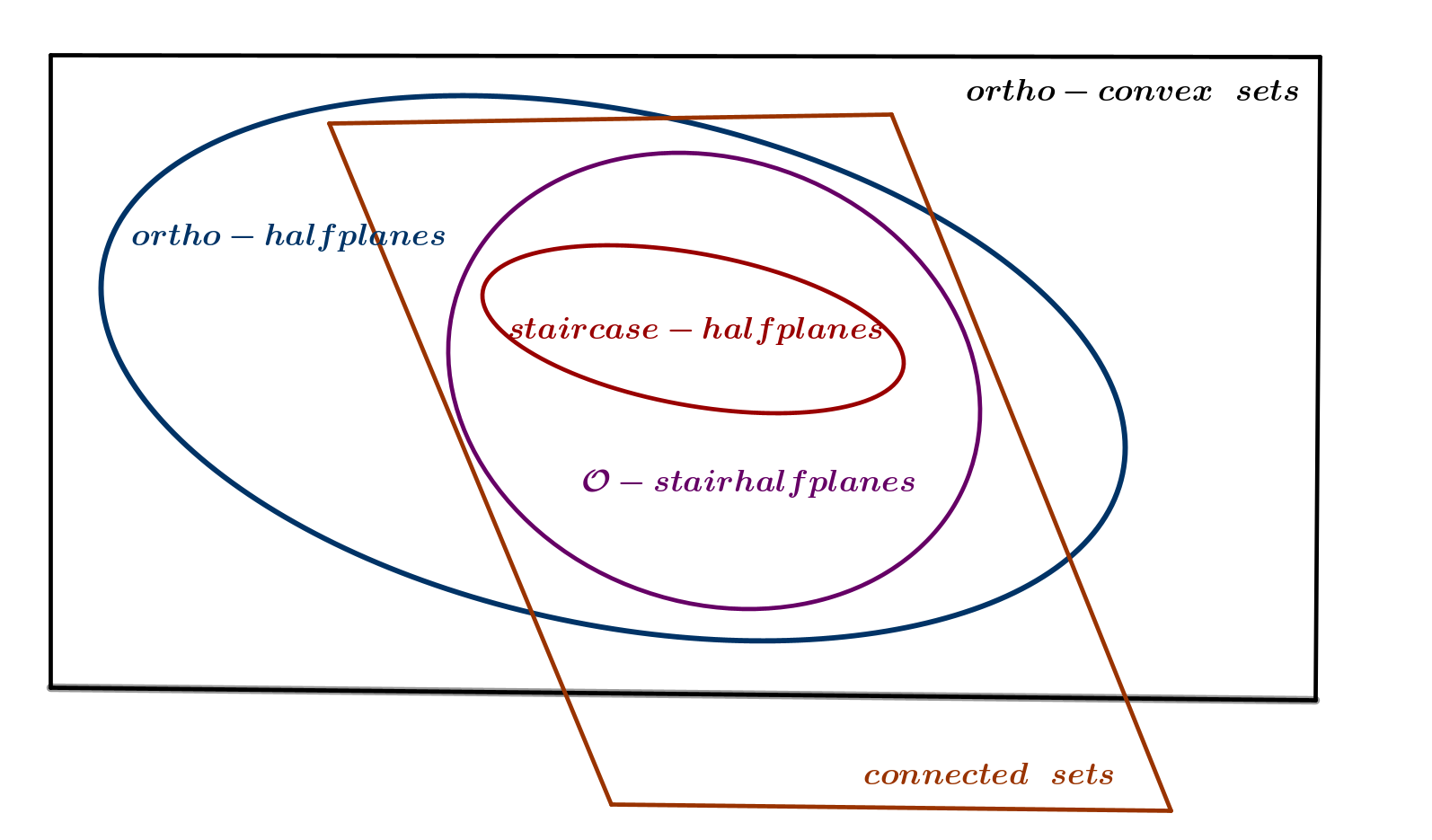}
		\caption{The relationship of three concepts of halpfplanes: staircase-halfplane, $\mathcal{O}$-stairhalfplane, and ortho-halfplane.}
		\label{fig:relationship-halfplane}
	\end{figure}

	For ortho-convex sets, there were two concepts of halfplanes  that are 	\textit{ortho-halfplane}~\cite{Fink1998} and $\mathcal{O}$\textit{-stairhalfplane}~\cite{RawlinsWood1991}. 
	An \textit{ortho-halfplane} is a closed set whose intersection with every axis-aligned line is empty, a ray or a line.
	Whereas an $\mathcal{O}$-\textit{stairhalfplane} is a region of the plane bounded by an ortho-convex path.
%
As shown in~\cite{Fink2012} and~\cite{RawlinsWood1991}, boundaries of $\mathcal{O}$-stairhalfplane can be  curved even though those of  ortho-halfplanes can be curved and disconnected. Whereas, the boundaries of staircase-halfplanes are  staircase lines that are polygonal lines.
	%
%
	%
	Since a staircase line is  ortho-convex, then every staircase-halfplane is an $\mathcal{O}$-stairhalfplane. Furthermore, by Lemma 5.8 in~\cite{Fink2012}, every $\mathcal{O}$-stairhalfplane is an ortho-halfplane. The relationship between three concepts of halpfplanes is shown in Fig.~\ref{fig:relationship-halfplane}.

\begin{definition}
	\label{def:separate-2sets}
	In $\mathbb{R}^2$, a staircase line \textit{separates} two sets $A$ and $B$ if one of two staircase-halfplanes  formed by the staircase line contains $A$ and the remaining one contains $B$.
	If neither $A$ nor $B$ does not intersect the staircase line, we say that the staircase line \textit{strictly separates}  $A$ and $B$.
\end{definition}

The ortho-convexity in $\mathbb{R}^n$ is generalized as follows
\begin{definition}[\cite{Fink2012}]
	\label{def:o-convex-Rn}
	A set $S \subset \mathbb{R}^n$ is said to be 
	\textit{ortho-convex}  if its intersection with any  axis-aligned line is connected.
\end{definition}

Clearly, the intersection of a family of ortho-convex sets is ortho-convex.
The intersection of all ortho-convex sets containing $S$ is called the \textit{orthogonally convex hull} of $S$, denoted by $\hbox{ortho-hull}(S)$.

Unlike traditional convex sets, an ortho-convex set is not necessarily path-connected, see Fig.~\ref{fig:examp-non-seperate}(i). Therefore the assumption of path-connectedness for ortho-convex sets is sometimes included in statements of the upcoming results. Furthermore, it is easy to show that combining the path-connectedness with the ortho-convexity in the plane leads to the simply connectedness.

\section {Blaschke-Type  Theorem for Ortho-Convex Sets in the Plane}
In this section, we introduce  Blaschke-type theorem for ortho-convex sets in the plane.
From now on, we denote the sub-path of a path $\gamma$ joining two points $u$ and $v$ of $\gamma$ by $\gamma(u,v)$.

\begin{lemma}
	\label{lem:sandwich-theorem}
	Let $\gamma$ be a path joining two points $a$ and $b$ in $\mathbb{R}^2$. Then $\gamma$ is ortho-convex if and only if it is $xy$-monotone.
	Moreover, if $\gamma$ is ortho-convex, we have
	\begin{align}
		\label{eq:sandwich-eq}
		\|a-b\| \le l(\gamma) \le \|a-b\|_1,
	\end{align}
	where $l(\gamma)$ is the length of $\gamma$ in Euclidean space, $\|.\|_1$ is $l^1$-norm and $\|.\|$ is Euclidean norm, i.e.,
	$\|a\|_1 = |a_x| +|a_y|, \, \|a\| = \sqrt{a_x^2+a_y^2}$, for $a=(a_x,a_y)$.
\end{lemma}
\begin{proof}
	($\Rightarrow$) Assume that $\gamma$ is ortho-convex and  not $xy$-monotone. W.l.o.g, assume that there are three points $u=(u_x,u_y), v=(v_x,v_y), w=(w_x,w_2) \in \gamma$ such that $u_x < v_x < w_x$, but $u_y <v_y$ and $v_y > w_y$.  Take $\alpha$ such that $\max\{u_y, w_y\} <\alpha <v_y $. Then the line $y = \alpha$ intersects with $\gamma$ in at least of two distinct points, say $m$ and $n$. By the ortho-convexness of $\gamma$ , we have $[m,n] \subset \gamma$. We take $\alpha$ changing between $\max\{u_y, w_y\}$ and $v_y$, we obtain another couple of distinct points $m'$ and $n'$ such that $[m',n'] \subset \gamma$. This contradicts the fact that $\gamma$ is a path. Therefore $\gamma$ is  $xy$-monotone.
	
	($\Leftarrow$) Suppose that $\gamma$ is $xy$-monotone and
	$l$ is a horizontal line such that 	$l \cap \gamma$ contains two distinct points $m=(m_x,m_y)$ and $n=(n_x,n_y)$, where $m_x <n_x$ (the case of vertical lines is similarly considered). Let $t = (t_x,t_y) \in \gamma(m,n)$ such that $t \neq m, t\neq n$. Then $m_x <t_x <n_x$. Thus $m_x -t_x<0$ and $n_x -t_x>0$.
	By the $xy$-monotonicity of $\gamma$, we have $(m_y-t_y)(n_y-t_y) \le 0.$ Note that $m_y =n_y$, it follows that $t_y=m_y=n_y$, then $t \in l$. Consequently $[m,n]=\gamma(m,n) \subset \gamma$. Hence $\gamma$ is ortho-convex.
	
	Now we are in a position to prove (\ref{eq:sandwich-eq}).
	Clearly, $l(\gamma) \ge \|a-b\|$. 
	Assume that $a_x \le b_x$ and  $\gamma$ is increasingly $xy$-monotone, i.e., if $a_x \le u_x \le v_x \le b_x$, then $u_y \le v_y$ (the case of decreasingly $xy$-monotonicity is similarly considered).
	We write $\gamma:[a_x,b_x]\to \mathbb{R}^2$.
	Then $l(\gamma) = \sup\limits_{\sigma} \sum\limits_{i=0}^{n-1} \|\gamma(t_i)-\gamma(t_{i+1})\|,$
	where $\sigma = \{t_i\}_{i=0}^{n}$ is a subdivision  of $[a_x,b_x]$. 
	Denote $\gamma^i = \gamma(t_i) =( \gamma^i_x, \gamma^i_y)$
	
	\begin{align*}
		\| & \gamma(t_i)-\gamma(t_{i+1})\|   = \sqrt{ (\gamma^i_x - \gamma^{i+1}_x)^2 + (\gamma^i_y - \gamma^{i+1}_y)^2} \\
		& = \sqrt{ \left( \gamma^i_x - \gamma^{i+1}_x + \gamma^i_y - \gamma^{i+1}_y \right) ^2 -2 (\gamma^i_x - \gamma^{i+1}_x )(\gamma^i_y - \gamma^{i+1}_y) } \\
		& \le \sqrt{ \left( \gamma^i_x - \gamma^{i+1}_x + \gamma^i_y - \gamma^{i+1}_y \right) ^2 }\\
		& =   \gamma^{i+1}_x - \gamma^i_x +  \gamma^{i+1}_y - \gamma^i_y \, \, \text{ as } \gamma \text{ is increasingly } xy \text{-monotone.}
	\end{align*}
	Consequently, 
	\begin{align*}
		l(\gamma) & \le \sup\limits_{\sigma} \sum\limits_{i=0}^{n-1} \left(\gamma^{i+1}_x - \gamma^i_x +  \gamma^{i+1}_y - \gamma^i_y \right) \\
		& = \sup\limits_{\sigma} \left(\gamma^n_x - \gamma^0_x + \gamma^n_y - \gamma^0_y\right) \\
		&= \gamma^n_x - \gamma^0_x + \gamma^n_y - \gamma^0_y \\
		& = \| \gamma^0 - \gamma^n \|_1 = \|a-b\|_1.
	\end{align*}
	Hence $\|a-b\| \le l(\gamma) \le \|a-b\|_1$. The proof is complete.
\end{proof}

\begin{lemma}
	\label{lem:o-convex-path-existance}
	Let $S$ be a closed, path-connected, and
	ortho-convex set in $\mathbb{R}^2$, $a,b \in S$. There exists a path joining $a$ and $b$ in $S$ which is ortho-convex. (Such a path is called an \textit{ortho-convex path}.)
\end{lemma}
\begin{proof}
		Due to the closeness and path-connectedness of $S$, there is the shortest path $\gamma$ joining $a$ and $b$ in $S$.
	Let $l$ be an axis-aligned line and $a_1,b_1 \in l \cap \gamma$.
	Since $\gamma(a_1,b_1)$ is shortest joining $a_1$ and $b_1$, we conclude that $[a_1,b_1] \subset S$, and therefore $\gamma$ is ortho-convex.
\end{proof}

\begin{lemma}
	\label{lem:sequence-length-paths}
	Suppose that $\{a_n\}$ and $\{b_n\}$ are sequences in $\mathbb{R}^2$, $a_n \rightarrow a$, $b_n \rightarrow b$ and $[a,b]$
	is an axis-aligned segment. 
	Let $\gamma_n$ be an ortho-convex path joining $a_n$ and $b_n$, for $n =1,2,\ldots$
	Then 
	$$l( \gamma_n)\rightarrow \|a-b\|  \, \text{ as } \, n \rightarrow \infty.$$
\end{lemma}

\begin{proof}
	By Lemma~\ref{lem:sandwich-theorem}, we have $\|a_n-b_n\| \le l(\gamma_n) \le \|a_n-b_n\|_1$. Then
	$\lim\limits_{n \to \infty} \|a_n-b_n\| \le \lim\limits_{n \to \infty} l(\gamma_n) \le \lim\limits_{n \to \infty} \|a_n-b_n\|_1$. Due to the continuity of $\|.\|$ and $\|.\|_1$, we have $\|a-b\| \le \lim\limits_{n \to \infty} l(\gamma_n) \le \|a-b\|_1$. Because $ \|a-b\| = \|a-b\|_1$, we obtain the required conclusion.
\end{proof}

\begin{lemma}
	\label{lem:Hausdorf-distance-paths}
	Suppose that $\{a_n\}$ and $\{b_n\}$ are sequences in $\mathbb{R}^2$, $a_n \rightarrow a$, $b_n \rightarrow b$ and $[a,b]$
	is an axis-aligned segment. 
	Let $\gamma_n$ ($\gamma$, resp.) be an ortho-convex path joining $a_n$ and $b_n$, for $n =1,2,\ldots$ ($a$ and $b$, resp.).
Then 
	$$d_\mathcal{H} \left( \gamma_n, \gamma \right) \rightarrow 0 \, \text{ as } \, n \rightarrow \infty.$$
\end{lemma}

\begin{proof} 
	Clearly, $\gamma =[a,b]$.
	 Suppose, contrary to our claim, that 
	$d_\mathcal{H} \left( \gamma_n, \gamma \right) \nrightarrow 0.$
	Then there is an $\epsilon > 0$ such that for each positive integer $k$ there exists an $n > k$ for
	which $d_\mathcal{H} \left( \gamma_n, \gamma \right) \ge 2 \epsilon$. Without loss of generality, we can assume that
	$d_\mathcal{H} \left( \gamma_n, \gamma \right) \ge 2 \epsilon$
	for all $n=1,2, \ldots.$
	 Since all convergent sequences are bounded,  there is  an axis-aligned rectangle, denoted by $R$, such that $\{a_n\}, \{b_n\} \subset R$ and $a,b \in R$. Obviously, all ortho-convex paths joining two points are closed and completely contained in an axis-aligned rectangle whose two opposite vertices are these two points. Then 
	 $\gamma$ and $\gamma_n$ are nonempty compact sets. There are $u_n \in \gamma_n$ and $v_n \in \gamma$ such that
	$$d(u_n, \gamma) = \sup_{y \in \gamma_n} d(y, \gamma) \, \text{ and } \, d(v_n, \gamma_n ) = \sup_{x \in \gamma}  d(x,\gamma_n ).$$
	Due to the definition of Hausdorff distance, either
	\begin{align}
		\label{eq:lem-limit1}
		d(u_n, \gamma) = d_{\mathcal{H}} ( \gamma_n, \gamma)
	\end{align}
	or
	\begin{align}
		\label{eq:lem-limit2}
		d(v_n, \gamma_n ) = d_{\mathcal{H}} ( \gamma_n, \gamma).
	\end{align}
	\textbf{Claim 1:} (\ref{eq:lem-limit1}) does not fulfill.
	Indeed, since $R$ is a compact set containing $\{u_n\}, \{v_n\}$, there are convergence subsequences $\{u_{n_j} \}$ and $\{v_{n_j} \}$. Let $u_0 = \lim\limits_{j \to \infty} u_{n_j}$  and $v_0 = \lim\limits_{j \to \infty} v_{n_j}$, then $u_0 \in R$ and $v_0 \in \gamma$.
	We have

	\begin{align*}
		\|a-u_0\|+\|u_0-b\| &= \lim\limits_{j \to \infty} \|a_{n_j}-u_{n_j}\|+ \lim\limits_{j \to \infty} \|u_{n_j}-b_{n_j}\| \text{ (by the continuity of } \|.\| \text{)}
		 \\
		  & \le  \lim\limits_{j \to \infty} l(\gamma_{n_j}(a_{n_j}, u_{n_j})) + \lim\limits_{j \to \infty}  l(\gamma_{n_j}(u_{n_j}, b_{n_j})) \text{ (due to Lemma~\ref{lem:sandwich-theorem})} \\ 
		& \le \lim\limits_{j \to \infty} \left[ l(\gamma_{n_j}(a_{n_j}, u_{n_j})) +  l(\gamma_{n_j}(u_{n_j}, b_{n_j})) \right] 
		\\
		& = \lim\limits_{j \to \infty}  l(\gamma_{n_j}(a_{n_j}, b_{n_j})) \text{ (since } u_{n_j} \text{ belongs to } \gamma_{n_j} \text{)}
		\\
		& \le  \|a-b \|_1    \text{ (due to Lemma~\ref{lem:sandwich-theorem})}  \\
		&= \|a-b\| \text{ (as } [a,b] \text{ is an axis-aligned segment)}
	\end{align*}
	Consequently, $u_0 \in [a,b] = \gamma$. Hence $d(u_{n_j},\gamma) \le \|u_{n_j} - u_0 \| \to 0$. 
	%
	Therefore (\ref{eq:lem-limit1}) does not hold for $n_j$ large enough. 
	
	\textbf{Claim 2:} (\ref{eq:lem-limit2}) does not fulfill. Assume, contrary to our claim, that for $n_j$ large enough, there is a natural number $j_1$ such that 
	\begin{align*}
		d(v_{n_j},\gamma_{n_j}) =  d_{\mathcal{H}} ( \gamma_{n_j}, \gamma) \ge 2 \epsilon, \, \, j \ge j_1.
	\end{align*}  
	It follows that 
	\begin{align*}
		\| v_{n_j} -y\| \ge 2 \epsilon \text{ for all } y \in \gamma_{n_j} \text{ and }
		j \ge j_1.
	\end{align*}
	Because $v_{n_j} \to v_0$, there exists a natural number $j_2 $ such that $\| v_{n_j} - v_0\| < \epsilon$ for $j \ge j_2$. Thus for  $j \ge \max \{j_1, j_2\}$ we obtain
	\begin{align}
		\label{eq:lem-limit3}
		\| v_0 -y \| \ge \| v_{n_j} -y \| - \| v_{n_j} - v_0 \| > 2 \epsilon - \epsilon = \epsilon \text{ for all } y \in \gamma_{n_j}.
	\end{align}
	More specifically, $\|v_0 - a_{n_j} \| > \epsilon$ and $\| v_0 - b_{n_j} \| > \epsilon$, for  $j \ge \max \{j_1, j_2\}$. As $a_{n_j} \rightarrow a$, $b_{n_j} \rightarrow b$, we get $v_0 \neq a$ and $v_0 \neq b$.
Therefore $0 <s:=\|a-v_0\| < \|a-b\|$. Furthermore, according to Lemma~\ref{lem:sequence-length-paths}, we get
	$l(\gamma_{n_j}) \rightarrow \|a- b\|$, there exists a natural number $j_3$ such that $l(\gamma_{n_j}) > s$ for $j > j_3$. Take $t_{n_j} \in \gamma_{n_j}$ satisfying $l(\gamma_{n_j}(a_{n_j}, t_{n_j})) = s$ for $j > j_3$. By the compactness of $R$, there is a convergent subsequence $\{t_{n_{j_k}}\}$ of $\{t_{n_j}\}$, $t_{n_{j_k}} \to t_0.$ We apply the argument above, with $u_0$ replaced by $t_0$, to obtain $t_0 \in [a,b]$. Using Lemma~\ref{lem:sequence-length-paths}, we have
	$$
	\|a-t_0\| = \lim\limits_{k \to \infty} l(\gamma_{n_{j_k}}(a_{n_{j_k}}, t_{n_{j_k}} )) = s.$$
	Thus $\|a-v_0\| =  \|a-t_0\|$, then $t_0 \equiv v_0$ and $ t_{n_{j_k}} \to v_0$ as $j \to \infty$. Note that $t_{n_{j_k}} \in \gamma_{n_{j_k}}$, this contradicts  (\ref{eq:lem-limit3}).
	
	Combining Claims 1 with 2, the proof is complete.
\end{proof}

\begin{corollary}
	\label{cor:closure-o-convex}
	\begin{itemize}
		\item[(a)] In $\mathbb{R}^2$, let $[a,b]$ be an axis-aligned segment and $c \in [a,b]$. Suppose that $\{a_n\}$, $\{b_n\}$ are two sequences satisfying $a_n \to a, b_n \to b$ and $\gamma_n$ is an ortho-convex path joining $a_n$ and $b_n$. Then there is $c_n \in \gamma_n, n=1,2,\ldots$ such that $c_n\to c$ as $n\to \infty$.
		\item[(b)] If $S$ is a path-connected and ortho-convex subset in $\mathbb{R}^2$, then so is  its closure ${\rm cl}S$.
	\end{itemize}
\end{corollary}
\begin{proof}
	(a) According to Lemma~\ref{lem:Hausdorf-distance-paths}, we have $d_\mathcal{H} \left( \gamma_n, \gamma \right) \rightarrow 0 \, \text{ as } \, n \rightarrow \infty.$
Since $\gamma_n$ is closed, there is $c_n \in \gamma_n$ such that $\|c-c_n\| = d(c,\gamma_n)$, for $n=1,2,\ldots$
Because $\|c-c_n\| \le d_\mathcal{H} \left( \gamma_n, \gamma \right)$, we get $c_n \to c$ as $n \rightarrow \infty.$
	
	(b)	Suppose that $l$ is an axis-aligned line and there are $a,b \in l \cap {\rm cl}S, a\neq b$. Then there are sequences $\{a_n\}, \{b_n\}$ in $S$ such that $a_n \to a$ and $b_n \to b$. By Lemma~\ref{lem:o-convex-path-existance}, there exists an ortho-convex path $\gamma_n$ in $S$ joining $a_n$ and $b_n$. 
	For $c \in [a,b]$, part (a) says that there is $c_n \in \gamma_n \subset S, n=1,2,\ldots$  such that $c_n \to c$ as $n \rightarrow \infty.$ Thus $c \in$ cl$S$. Hence $[a,b] \subset$ cl$S$, and cl$S$ is ortho-convex.
\end{proof}

\begin{lemma}
	\label{lem:limit-sequence-o-convex}
	Let $\{S_i\}$ be a family of nonempty compact, path-connected, and ortho-convex subsets in $\mathbb{R}^2$
	and suppose that $\{S_i\}$ converges to a nonempty compact set $S$, i.e., $d_{\mathcal{H}} ( \mathcal{S}_n, \mathcal{S}) \to 0$. Then $S$ is also ortho-convex.
\end{lemma}

\begin{proof}
Suppose that $l$ is an axis-aligned line and there are $a,b \in l \cap S, a\neq b$. Since $S_n$ is closed, there is $a_n \in S_n$ such that $\|a-a_n\| = d(a,S_n)$, for $n =1,2,\ldots$ Because $d(a,S_n) \le \sup_{x \in S}d(x,S_n) \le d_{\mathcal{H}} ( \mathcal{S}_n, \mathcal{S}) \to 0$, we have $a_n \to a$. Similarly there is a sequence $\{b_n\}$ such that $b_n \in S_n, b_n \to b$.
Using Lemma~\ref{lem:o-convex-path-existance}, we get an ortho-convex path $\gamma_n$ in $S_n$ joining $a_n$ and $b_n$, for $n=1,2,\ldots$
By Lemma~\ref{lem:Hausdorf-distance-paths}, we obtain $d_{\mathcal{H}} ( \gamma_n, \gamma) \to 0$. For $c \in [a,b]$, according to Corollary~\ref{cor:closure-o-convex} part (a), there exists $c_n \in \gamma_n \subset S_n$ such that $c_n \to c$. 
As $S$ is closed, we take a point $d_n \in S$ such that $\|c_n -d_n\| =d(c_n,S)$, for $n =1,2,\ldots..$ Then $\|c_n -d_n\| \le d_{\mathcal{H}} ( \mathcal{S}_n, \mathcal{S})$ and thus,
$$\|d_n -c\| \le \|d_n -c_n\|+\|c_n-c\| \le d_{\mathcal{H}} ( \mathcal{S}_n, \mathcal{S}) +\|c_n -c\| \to 0.$$
	Because of the closeness of $S$, it follows that $c \in S$. Therefore $[a,b] \subset S$ and $S$ is ortho-convex.
%
%
\end{proof}

 It should be noticed that if the assumption
of closeness in the above lemma is dropped, then we consider $\{{\rm cl}S_n\}$ and ${\rm cl}S$ and
conclude that ${\rm cl}S$ is ortho-convex.
Next, we prove the main result of this section.

\begin{theorem}[Blaschke-type theorem]
	\label{theo:Blaschke}
	Let $\mathcal{F}$ be a 
	uniformly bounded infinite collection of nonempty compact, path-connected, and ortho-convex sets in $\mathbb{R}^2$. Then $\mathcal{F}$ contains a subsequence that converges to a nonempty compact and ortho-convex set.
\end{theorem}

\begin{proof}
	Let $\mathcal{K}$ be the collection of all nonempty compact sets in $\mathbb{R}^2$. Then $\left( \mathcal{K},d_{\mathcal{H}}\right)$ is a metric space. Since $\left(\mathbb{R}^2, \|.\| \right)$ is a complete metric space, so is $\left( \mathcal{K},d_{\mathcal{H}}\right)$, due to the theorem of~\cite{Price1940}, page 1. Because $\mathcal{F}$ is a  
	uniformly bounded infinite collection of elements in $\left( \mathcal{K},d_{\mathcal{H}}\right)$, then $\mathcal{F}$ contains a subsequence that converges to a nonempty compact set $S$ in $\left( \mathcal{K},d_{\mathcal{H}}\right)$. By Lemma~\ref{lem:limit-sequence-o-convex}, $S$ is ortho-convex. Hence the proof is complete.
	\end{proof}
	
\section {Separation of Disjoint Closed, Path-connected, and Ortho-Convex Sets in the Plane}

Some of the most important applications of convex sets involve the problem of
separating two convex sets by a hyperplane. 
Herein, a staircase-halfplane plays the same role as a traditional halfplane in $\mathbb{R}^2$ do for convex sets. 




\begin{lemma}
	\label{lem:separate-pnt-set}
	Let $S$ be  closed, path-connected,  and ortho-convex. If $p \notin S$, then there exists  a staircase line strictly separating $S$ and $p$.
\end{lemma}

\begin{proof}
	
	If $S = \emptyset$, the proof is trivial. Assume that $S$  is non-empty.
	Let $d = d(p,S)$.
	Because $S$ is closed and $p \notin S$, then  $d >0$, and there exists $a \in S$ such that $d = \| p-a \|$.
	We construct a grid whose size is $\frac{d}{2\sqrt{2}}$.
	Let $P$  be the union of all grid cells containing $p$.
	Then $P$ is a box including one, two,  or $2\times 2$ grid cells.

	\textbf{Claim 1.} $S \cap P = \emptyset$. 	
	Indeed, suppose, contrary to our claim, that there exists $ b \in S \cap P$. Then 
	\begin{align}
		\label{eq:1}
		d = d(p,S) \le \| p-b\|
	\end{align}
	By the construction $P$ and
	$b \in P$, then $\| p-b\|$ does not exceed the diameter of one grid cell.
	Therefore 
	\begin{align}
		\label{eq:2}
		\| p-b\| \le \sqrt{2}.\frac{d}{2\sqrt{2}} = \frac{d}{2} <d
	\end{align}
	Combining~(\ref{eq:1}) and ~(\ref{eq:2}), we get a contradiction and thus the claim is proved.
	
	Let $c,d,e$, and $f$ be four vertices of $P$. Corresponding to each vertex, there is one staircase line that passes through two adjacent edges of $P$. Staircase lines determine four  quadrants\footnote{See~\cite{AnHuyenLe2021} for the definition of a quadrant.} together with the staircase lines containing the entire $P$, see Fig~\ref{fig:4quadrants}(ii).
	
	\begin{figure}[htp]
	\begin{subfigure}[b]{0.5\textwidth}
		\centering
		\includegraphics[scale=0.25]{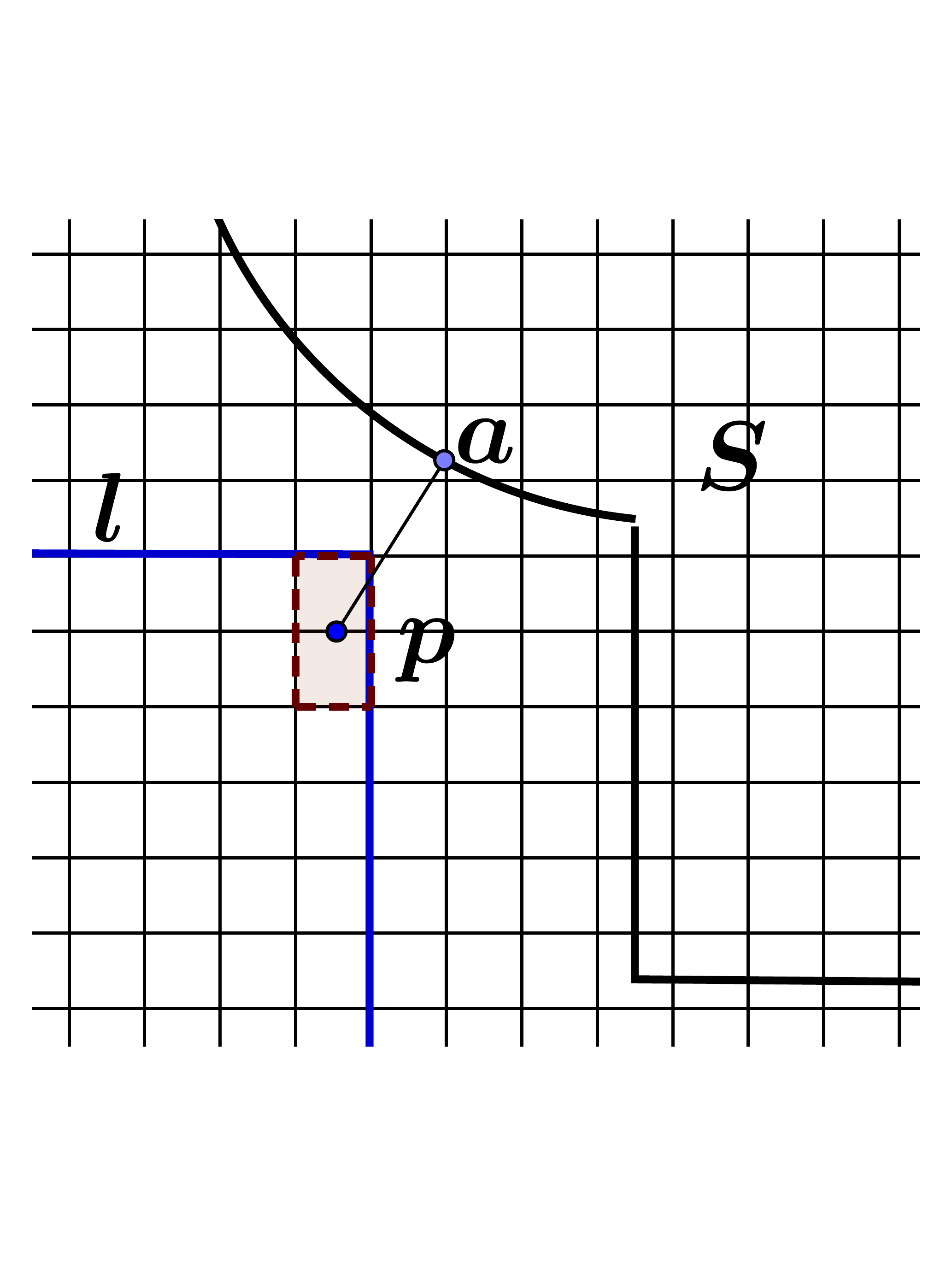} 
		\caption*{(i)}
	\end{subfigure}%
	\begin{subfigure}[b]{0.5\textwidth}
		\centering
		\includegraphics[scale=0.3]{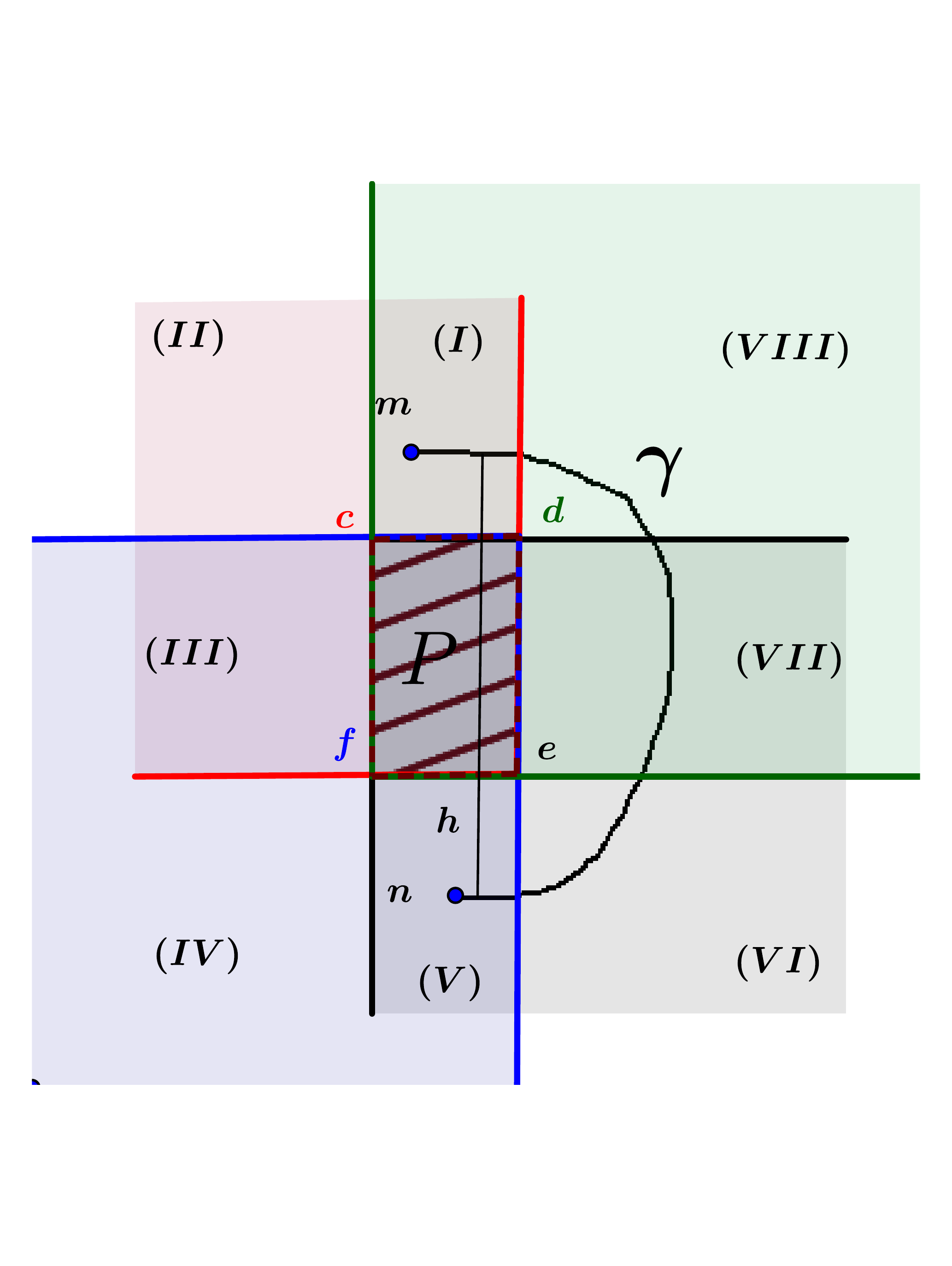}
		\caption*{(ii)}
	\end{subfigure}
	\caption{Illustration of the proof Lemma~\ref{lem:separate-pnt-set}.}
	\label{fig:4quadrants}
			\end{figure}  
	

	\textbf{Claim 2.} There is at least one quadrant intersecting $S$ in the empty set.
	
	Indeed, assume that all quadrants intersect $S$ in  nonempty sets. 
	We draw lines passing through the edges of $P$ and dividing the plane into eight
	closed regions, denoted by (I)-(VIII), see Fig~\ref{fig:4quadrants}(ii). 
	We prove that 
	\begin{align}
		\label{eq:5}
		\text{If }S \cap (I) \neq \emptyset \text{ then }S \cap (V) = \emptyset. 
	\end{align}
	
	W.l.o.g. suppose that (\ref{eq:5}) is not true, then there exist $m  \in S \cap (I)$ and $n \in S \cap (V)$.
	Because $S$ is path-connected, there is a path $\gamma$ joining $m$ and $n$ such that $\gamma \subset S$. Then there is a horizontal line segment $h$ joining two points of $\gamma$ such that $h \cap P \neq \emptyset$. Since $S$ is ortho-convex, $h \subset S$. Therefore $S \cap P \neq \emptyset$, which contradicts the proven result in Claim 1. 
	
	Because of the equality of (I), (III), (V), and (VII) regions, (\ref{eq:5}) holds for the remaining regions.
	Therefore we can assume that $S\cap \left( (III) \cup (V) \right) = \emptyset$. By the path-connectedness of $S$, we obtain either 
	$B \subset (IV) $ or $B \cap ((III) \cup (IV) \cup (V)) = \emptyset$.
	If  $S \subset (IV)$, 
	the staircase line corresponding to the quadrant having the vertex at $f$ strictly separates $S$ and $p$. 
	If $B \cap ((III) \cup (IV) \cup (V)) = \emptyset$, the staircase line corresponding to the quadrant having the vertex at $d$ strictly separates $S$ and $p$. 
	The proof is complete.
\end{proof}

\begin{figure}[htp]
	\begin{subfigure}[b]{0.45\textwidth}
		\centering
		\includegraphics[scale=0.25]{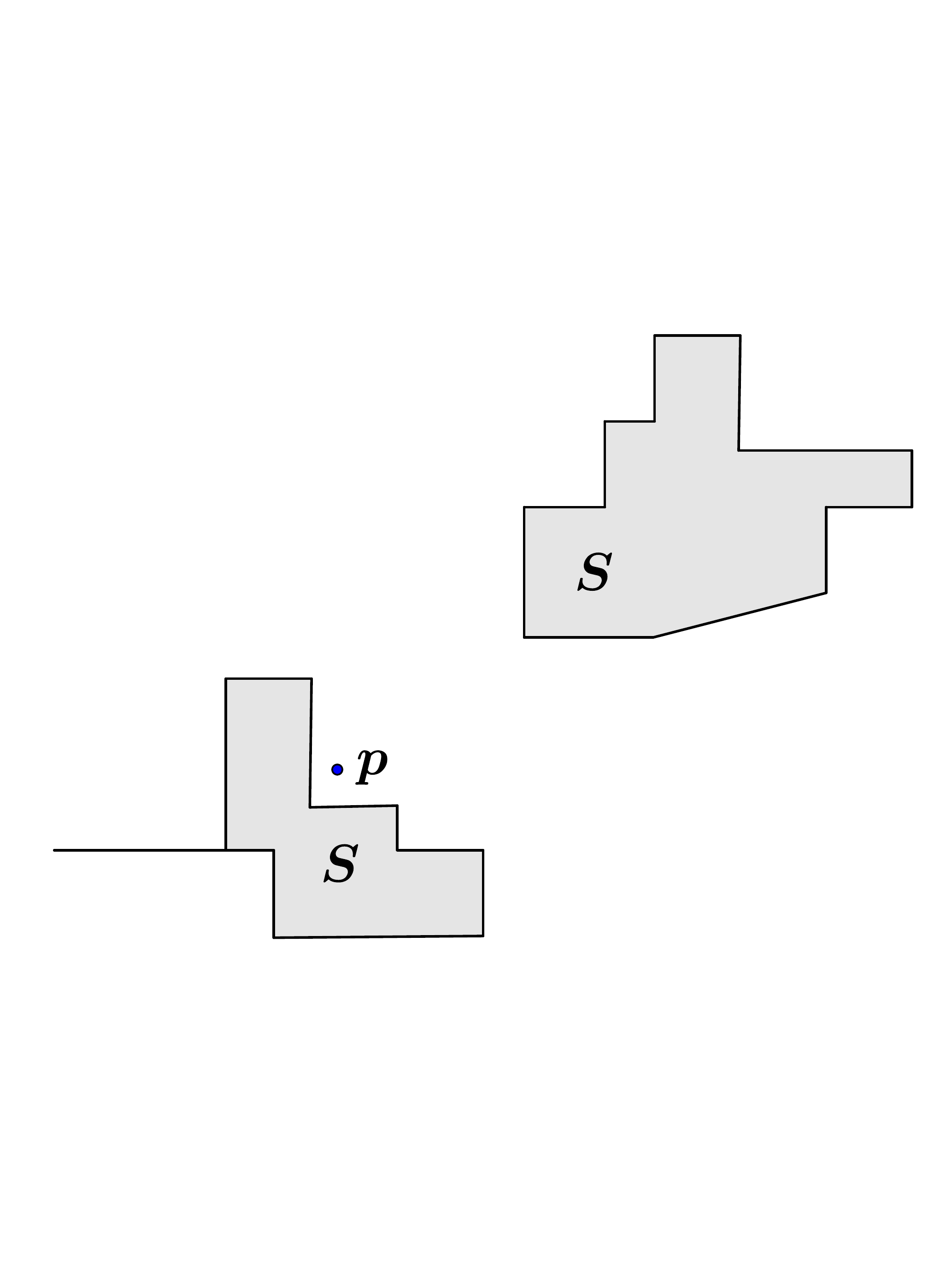}
		\caption*{(i)}
	\end{subfigure}%
	\begin{subfigure}[b]{0.45\textwidth}
	\centering
	\includegraphics[scale=0.2]{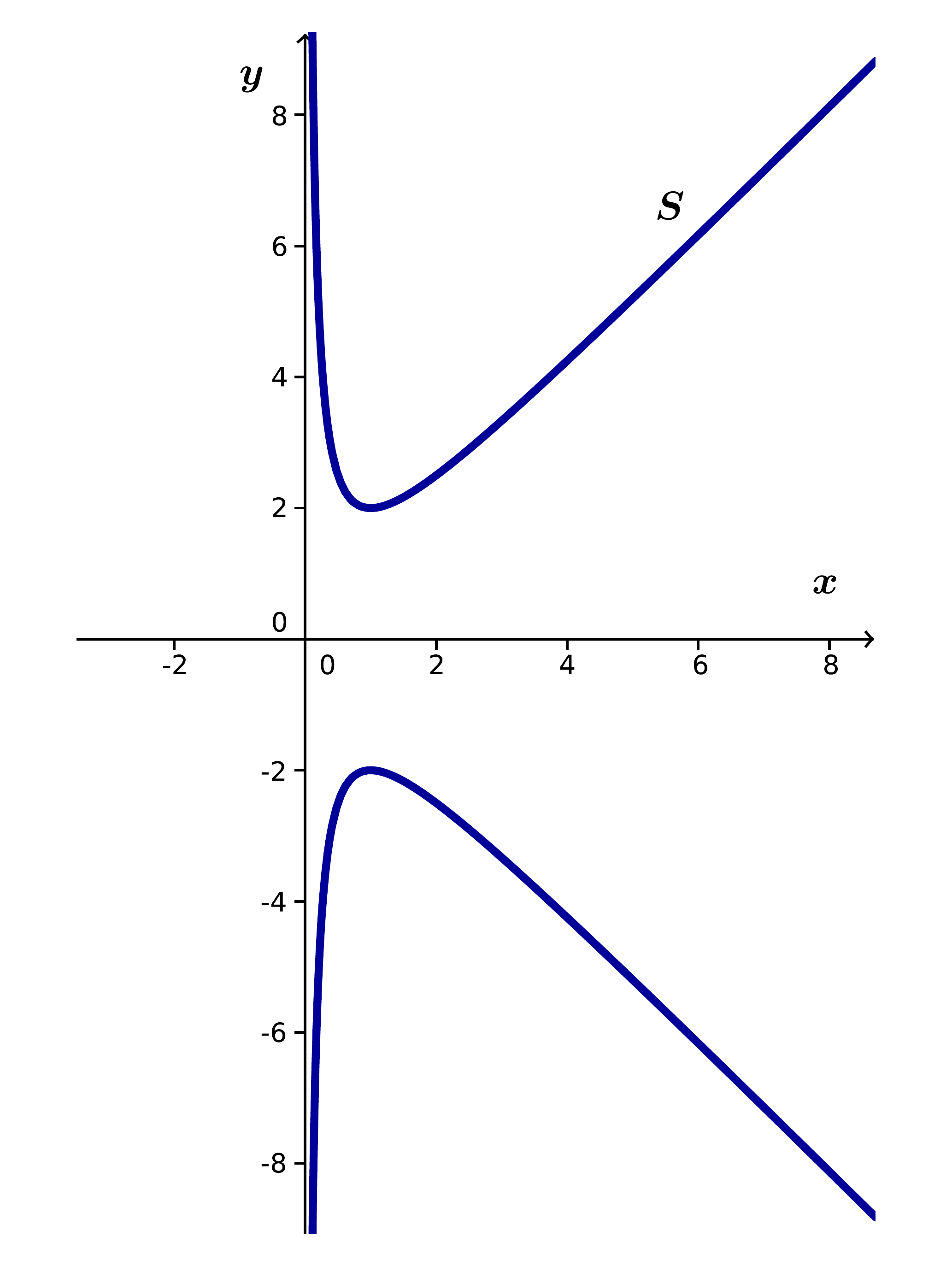}
	\caption*{(ii)}
\end{subfigure}
	\caption{(i):  $S$ is ortho-convex but not path-connected, and there is no staircase line separating $p$ and $S$ and (ii): the orthogonally convex hull of a closed set $S$ is not closed.}
	\label{fig:examp-non-seperate}
\end{figure}

Lemma~\ref{lem:separate-pnt-set} does not hold for sets that are not path-connected. Fig~\ref{fig:examp-non-seperate}(i) shows that 
there is no staircase line  separating $S$ and  $p$, although $S$ is closed and ortho-convex, $p \notin S$.



\begin{theorem}
	\label{theo:separate-2sets}
Suppose that $A$ and $B$ be two disjoint compact, path-connected, and ortho-convex sets. Then there exists  a staircase line strictly separating $A$ and $B$.
\end{theorem}
\begin{proof}
	If $A$ or $B$ is empty, the proof is trivial. Assume that $A$ and $B$ are non-empty.
	Let $d = d(A,B)$.
	Since $A$ and $B$ are two compact sets, $d >0$.
	We construct a grid whose size is $\frac{d}{2\sqrt{2}}$.
	Let $A^g$  be the union of all grid cells containing any element of $A$.
	As $A$ is compact, path-connected, and ortho-convex,  $A^g$ is clearly a simple polygon which is ortho-convex.
	
	We next show that $B \cap A^g = \emptyset$.
	Indeed, suppose, contrary to our claim, that there exists $ u \in B \cap A^g$. Then 
	\begin{align}
		\label{eq:3}
		d = d(A,B) \le d(u,A)
	\end{align}
	By the construction $A^g$ and
	$u \in A^g$, then $d(u,A)$ does not exceed the diameter of one cell. Thus 
	\begin{align}
		\label{eq:4}
		d(u,A) \le \sqrt{2}.\frac{d}{2\sqrt{2}} = \frac{d}{2} < d
	\end{align}
	Combining~(\ref{eq:3}) and ~(\ref{eq:4}), we get a contradiction. Thus $B \cap A^g = \emptyset$.

	Let $P$ be the minimal axis-aligned rectangle such that $P$ contains $A^g$. Let $a,b,c$ and $d$ be four vertices of $P$. 
	Similarly to the proof of Lemma~\ref{lem:separate-pnt-set} (Claim 2), we draw lines passing through edges of $P$ and dividing the plane into eight
	closed regions, denoted by (I)-(VIII), see Fig 5.
	
	By the same way, we can 
	show that if $B \cap (I) \neq \emptyset$  then $B \cap (V) = \emptyset$. Because of the equality of (I), (III), (V), and (VII) regions, the preceding argument can be applied to the remaining regions.
	Therefore we can assume that
	\begin{align}
		\label{eq:6}
		B\cap \left( (III) \cup (V) \right) = \emptyset. 
	\end{align}
	
	We take the points
	$m,n,p$, and $q$ belonging to $P$ such that they are respectively the highest, lowest  leftmost and rightmost of $P$. Let $P_1, P_2, P_3$ and $P_4$ be the sub-rectangles having three vertices $m, b$, and $q$; $q,c$, and $n$; $n,d$, and $p$; $p,a$, and $m$, respectively. They  can degenerate to a point or a line segment.

	If $B \cap P = \emptyset$, as in the proof of Lemma~\ref{lem:separate-pnt-set}, either 
	$B \subset (IV) $ or $B \cap ((III) \cup (IV) \cup (V)) = \emptyset$.
	Then  the staircase line corresponding to the quadrant having the vertex at either $d$  or $b$ strictly separates $S$ and $p$. 
	Otherwise, $B \cap P \neq \emptyset$. Then  $B$ intersects  only one of four rectangles $P_i, (1 \le i \le 4)$ by the path-connectedness and ortho-convexity of $B$.
	By (\ref{eq:6}) and the path-connectedness of $B$,
	assume, w.l.o.g, that $B \cap P_1 \neq \emptyset$. Then $P_1$ is a non-degenerate rectangle.
	\begin{figure}
		\centering
		\includegraphics[width=0.58\linewidth]{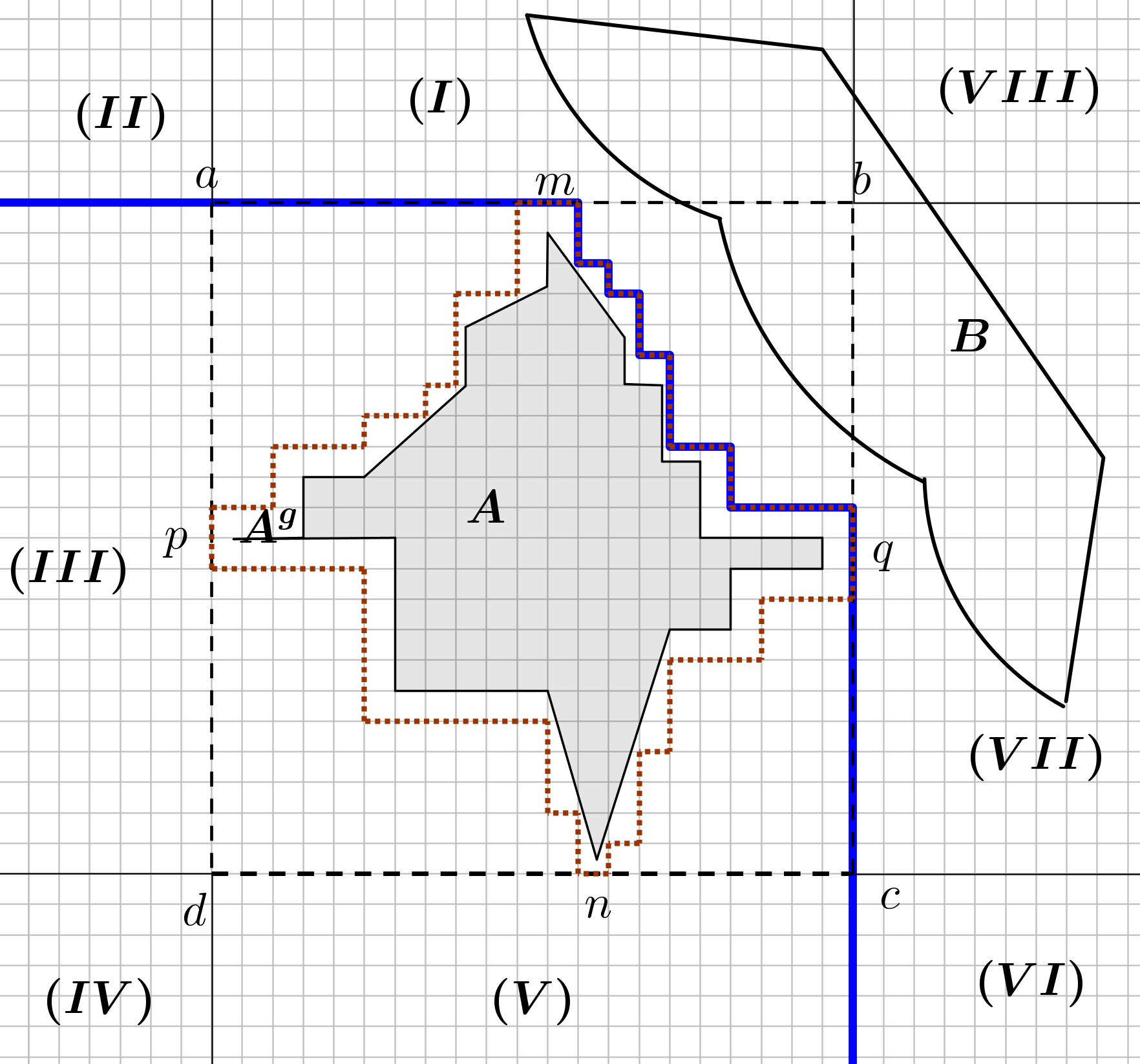}
		\caption{Illustration of the proof Theorem~\ref{theo:separate-2sets}.}
		\label{fig:seperate-2sets}
	\end{figure}
	Let $l$ is the path formed by the ray $ma$, the boundary of $P$ from $m$ to $q$, and the ray $qc$. (If $m=a$ or $q=c$, then $ma$ and $qc$ are understood as the opposite rays of $mb$ and $qb$). By the ortho-convexity of $P$, then $l$ is a staircase line. Moreover, $l$ strictly separates $A$ and $B$. The proof is complete.
\end{proof}

%

Since a staircase line is the specified case of an ortho-convex path, the separation of two disjoint ortho-convex sets by a staircase line leads to the separation of two disjoint ortho-convex sets by an ortho-convex path as shown in~\cite{Dulliev2017}. Note that the concept of an ortho-convex path is equivalently that of a monotone curve in~\cite{Dulliev2017}.

\section{Representing Closed Ortho-convex Sets by the Intersections of Staircase-Halfplanes}

Fink and Wood~\cite{Fink1998} represented a connected, closed, and ortho-convex sets by the intersections of ortho-halfplanes
as follows
\begin{lemma}[Lemma 2.3 in~\cite{Fink1998}]
	A connected, closed set $S$ is ortho-convex if and only if it is the intersection of ortho-halfplanes.
\end{lemma}

For staircase-halfplanes case, we obtain the following

\begin{theorem}
	\label{theo:represent-close-o-conv}
	A path-connected set is closed and ortho-convex if and only if it is the intersection of a family of staircase-halfplanes.
\end{theorem}
\begin{proof}
	$(\Leftarrow)$	Suppose that a set $S$ is the intersection of a family of staircase-halfplanes. Because the intersection of ortho-convex sets is ortho-convex and 	every
	staircase-halfplanes is ortho-convex, so is  their intersection. Furthermore, the intersection of a family of closed sets is closed, then $S$ is closed.
	
	$(\Rightarrow)$	Suppose that  $S$ is a closed, path-connected,  and ortho-convex set. For all $p \notin S$, according to Lemma~\ref{lem:separate-pnt-set}, there exists a staircase line strictly separating $S$ and $p$.
	Let $\mathcal{F}$ be the family of  staircase-halfplanes containing $S$ corresponding to each $p \notin S$.
	We prove that $S = \bigcap_{L \in \mathcal{F}} L$.
	Clearly $S \subset \bigcap_{L \in \mathcal{F}} L$.
	Take $a \notin S$, we show that  $a \notin \bigcap_{L \in \mathcal{F}} L$. Indeed, by Lemma~\ref{lem:separate-pnt-set},
	there is a staircase-halfplane $L \in \mathcal{F}$ containing $S$ but not $a$. Then $a \notin \bigcap_{L \in \mathcal{F}} L$. Thus $ \bigcap_{L \in \mathcal{F}} L \subset S$. Hence $S = \bigcap_{L \in \mathcal{F}} L$.	
\end{proof}

It should be noticed that if the assumption
of path-connectedness and closeness in the above theorem is dropped, then we do not obtain the required conclusion (see Example~\ref{exam_not_path-connected_not_closed}).

\begin{example}
\label{exam_not_path-connected_not_closed}
	%
	Example 20, page 142 in~\cite{Gelbaum2003} showed a set $A$ that is dense in $[0,1] \times [0,1]$, and the intersection of $A$ and any axis-aligned line has at most one point. Clearly $A$ is ortho-convex but $A$ is neither  path-connected nor closed. 
	Since all staircase-halfplanes are closed, the  intersection of any family of  staircase-halfplanes containing $A$ is 
	closed.
	Therefore the intersection of any family of  staircase-halfplanes containing $A$ is not $A$ as required.
\end{example}

	\section{Some Topological Properties of Ortho-Convex Sets in $\mathbb{R}^n$}\label{DefsAndProps}

In $\mathbb{R}^n$, let $\mathcal{H}_i=\{x=(x_k)_{k=1}^n \in \mathbb{R}^n|\,x_i = 0\}$, for $i=1,2,\ldots n$.

\begin{lemma}
	\label{lem:formula-o-convex}
	Let $S \subset \mathbb{R}^n$. The following three conditions are equivalent
	\begin{itemize}
		\item[(a)] $S$ is ortho-convex;
		\item[(b)]  The intersection of $S$ and each   hyperplane which is parallel to  $\mathcal{H}_i$  is ortho-convex, $i=1,2,\ldots n$.
		\item[(c)] For each $i=1,2,\ldots n$, we take two arbitrary points $a,b \in S, a = (a_1,a_2, \ldots, a_n), b = (b_1, b_2, \ldots, b_n)$ satisfying 
		
		\begin{align}
			\label{eq:convex-condition}
			\mathop  \sum_{\substack{1 \le k \le n \\ k \neq i}}(a_k -b_k)^2 =0,
		\end{align}  
		then  $\lambda a +(1- \lambda)b \in S$, for all  $\lambda \in [0,1]$.
	\end{itemize} 
	
\end{lemma}
\begin{proof}
	(a)$\Leftrightarrow$(b) is true, according to Theorem4.6 in~\cite{Fink2012}.
	\\
	(a)$\Leftrightarrow$(c) is obvious since two points $a,b \in S \cap l$ if and only if $a,b \in S$ and their coordinates satisfy (\ref{eq:convex-condition}), where $l$ is an arbitrary axis-aligned line.
\end{proof}


	\begin{lemma}\footnote{Incidentally,  this result was published in~\cite{Dulliev2017}. That was informed by the author of~\cite{Dulliev2017} after our submission.}
		\label{lem:interior-convex}
		The interior of an ortho-convex set is ortho-convex.
	\end{lemma}

	\begin{proof}
		Let $S$ be an  ortho-convex set and $l$ be an axis-aligned line. For two distinct points $a,b \in l \cap \text{int}S$, we prove $[a,b] \subset \text{int}S$. 
		 Due to $a,b \in \text{int}S$, there are two  balls $B(a,r)\subset S$ and $B(b,r) \subset S$, where $r >0$, see Fig.~\ref{fig:join-two-balls}. Take any $p \in [a,b]$, $p = (p_1,p_2, \ldots, p_n)$. W.l.o.g.  suppose that $l$ is parallel to the axis spanned by $e_i=(\underbrace{0,...,0,1}_{i},0\ldots,0)$.
		\begin{figure}[htp]
			\centering
			\includegraphics[width=0.42\linewidth]{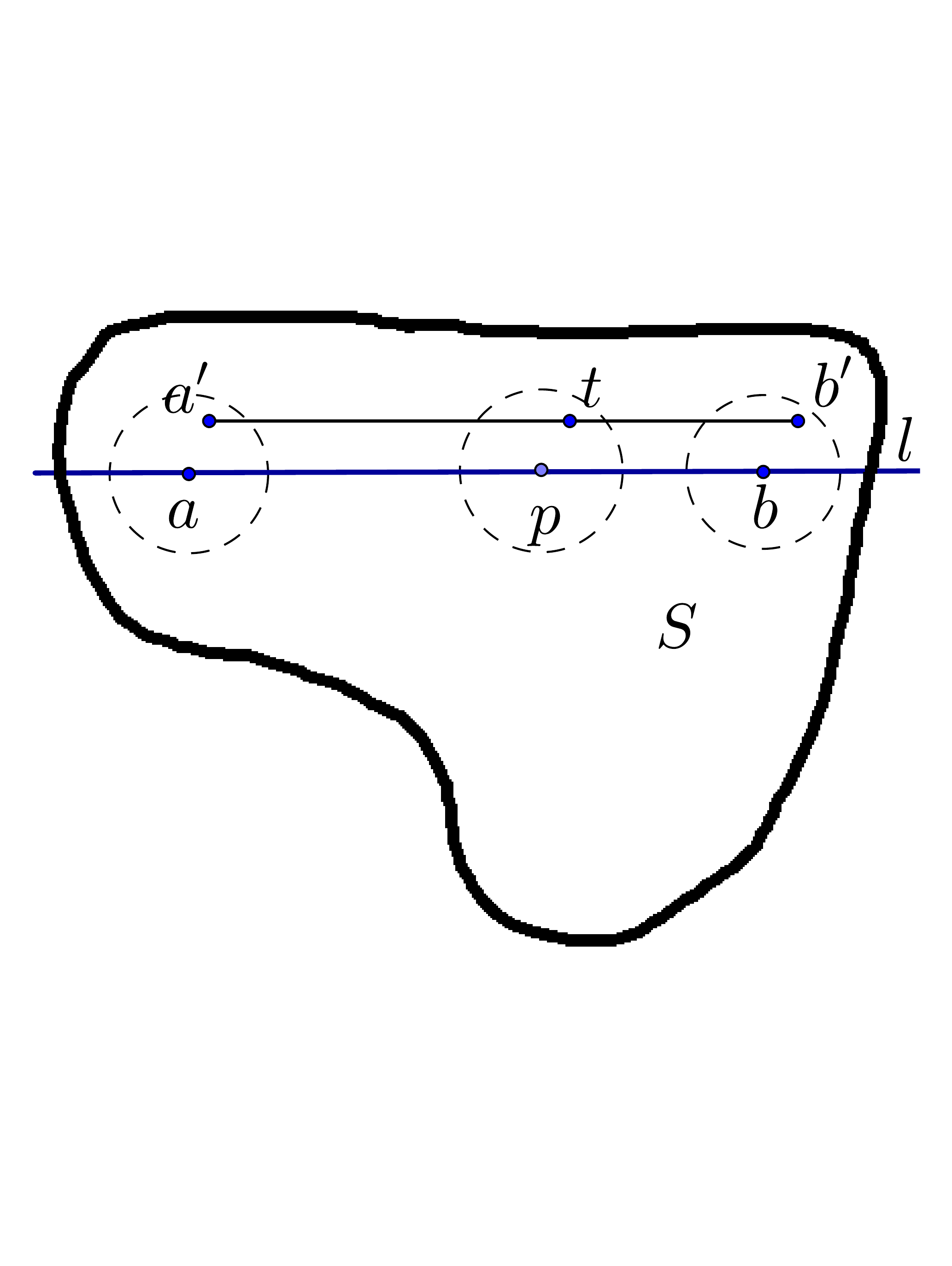}
			\caption{Illustration of the proof of Lemma~\ref{lem:interior-convex}.}
			\label{fig:join-two-balls}
		\end{figure}
		 Then $a = (p_1,p_2, \ldots,p_{i-1},\alpha,p_{i+1}, \ldots, p_n) $ and \break $b = (p_1,p_2, \ldots,p_{i-1},\beta,p_{i+1}, \ldots, p_n) $. By $p \in [a,b]$, there is $\lambda_0 \in [0,1]$ such that
		\begin{align}
			\label{eq:int-open1}
			p_i=\lambda_0 \alpha + (1-\lambda_0)\beta,
		\end{align}
		For all $t =(t_1,t_2, \ldots, t_n) \in B(p,r)$, let $a'=(t_1,t_2, \ldots,t_{i-1},\alpha,t_{i+1}, \ldots, t_n)$ and $b'=(t_1,t_2, \ldots,t_{i-1},\beta,t_{i+1}, \ldots, t_n)$.
		Due to $t  \in B(p,r)$, we have
		\begin{align}
			\label{eq:int-open2}
			\|p-t \| = \sum_{k=1}^{n} (p_k -t_k)^2 <r.
		\end{align} 
		We have 
		\begin{align}
			\label{eq:int-open3}
			\|a-a'\|= & (p_1 -t_1)^2 + \ldots + (p_{i-1} -t_{i-1})^2 + (\alpha -\alpha)^2+ (p_{i+1} -t_{i+1})^2 \\ \notag
			&+ \ldots + (p_n -t_n)^2 \\ \notag
			\le &\sum_{k=1}^{n} (p_k -t_k)^2 <r, \text{ due to (\ref{eq:int-open2}).}
		\end{align}
		Therefore $a' \in B(a,r) \subset S$. Similarly, we have $b' \in B(b,r)\subset S$.
		On the other hand, according to (\ref{eq:int-open1}), we have $t =\lambda_0 a' + (1-\lambda_0)b'$. 
		Clearly, coordinates of $a'$ and $b'$ satisfy (\ref{eq:convex-condition}).
		Because of the ortho-convexity of $S$ and Lemma~\ref{lem:formula-o-convex}(c), we get $t \in S$. Hence $B(p,r) \subset S$,  $p \in \text{int}S$ and thus $[a,b] \subset \text{int}S$. The proof is complete.
	\end{proof}


\begin{theorem}
	\label{theo:o-convex-open-set}
	The orthogonally convex hull of an open set is open.
\end{theorem}
\begin{proof}
	Let $S$ be   an open set and $A = \hbox{ortho-hull}(S)$. Then $S \subset A$.
	Since $S$ is open, we have $S \subset \text{int}A$. Furthermore, by Lemma~\ref{lem:interior-convex}, $\text{int} A$ is ortho-convex. Therefore $\text{int} A$ is ortho-convex set containing $S$. Because $A$ is the intersection of all ortho-convex sets containing $S$, then $A \subset \text{int}A$. On the other hand, $A = \text{int}A$, and thus $A$ is open. The proof is complete.
\end{proof}

Lemma~\ref{lem:interior-convex} and Theorem~\ref{theo:o-convex-open-set} do not hold for closeness. It means that  if $S$ is ortho-convex, then ${\rm cl}S$ may not be ortho-convex (see Example~\ref{exam_closure_not_o-convex}) and if $S$ is closed, then $\hbox{ortho-hull}(S)$ may not be closed (see Example~\ref{exam_not_close}).
However, if $S$ is a path-connected and ortho-convex set in the plane, then so is its closure cl$S$, see Corollary~\ref{cor:closure-o-convex} (b).

\begin{example}
	\label{exam_closure_not_o-convex}
	Let	$S = \{(x,0) \in \mathbb{R}^2 |\, -1<x <0 \} \cup \{(x,1) \in \mathbb{R}^2 |\, 0<x <1 \} $.
	Clearly
	$S$ is ortho-convex. Whereas ${\rm cl}S = \{(x,0) \in \mathbb{R}^2 |\, -1 \le x \le 0 \} \cup \{(x,1) \in \mathbb{R}^2 |\, 0 \le x \le 1 \} $ is not ortho-convex, since the line $x=0$ intersects $S$ in a set of two disjoint points $\{(0,0),(0,1)\}$.
	
\end{example}

\begin{example}
	\label{exam_not_close}
	Let	$S = \{(x,y) \in \mathbb{R}^2 |\, y^2= (\frac{x^2+1}{x})^2 , x>0  \}$, see~Fig.~\ref{fig:examp-non-seperate}(ii).
	Clearly
	$S$ is closed, but $\hbox{ortho-hull}(S) = \{ (x,y) \in \mathbb{R}  | \, x >0 \}$ is open.
	Indeed, let $A = \{ (x,y) \in \mathbb{R}^2  | \, x >0 \}$.
	For all $p \in A$, if $p \in S$, then $p \in \hbox{ortho-hull}(S)$. If $p \notin S$, there exists an axis-aligned line intersecting $S$ in a line segment containing $p$, and thus $p \in \hbox{ortho-hull}(S)$. Therefore $ A \subset \hbox{ortho-hull}(S)$. Conversely, for all $p \notin A$, any vertical or horizontal line either  intersects $S$ in a line segment not containing $p$, or does not. Then $ \hbox{ortho-hull}(S) \subset A$. Hence $ A = \hbox{ortho-hull}(S)$. 
	
\end{example}



\section{Concluding Remarks}

	We have provided some similar results as in convex analysis, that is Blaschke-type theorem and the separation  for ortho-convex sets in the plane. Therefore a closed and ortho-convex set is represented by the intersection of staircase-halfplanes. 
	We then introduce some topological properties of ortho-convex sets in 	$\mathbb{R}^n$. However, in $\mathbb{R}^n (n>2)$, some  properties of ortho-convexity do not remain such as  an ortho-convex path 
	 may not be monotone with respect to directions of coordinate axes, or a path-connected and ortho-convex set may be not simply connected.
There arises a question  ``Does  Blaschke-type theorem and the separation property
 for ortho-convex sets hold true in $\mathbb{R}^n (n >2)$?" This will be the subject of another paper.
In addition, because ortho-convexity is a special case of restricted orientation convexity which are discussed in~\cite{Alegria2022,Fink2012}, the results presented in this paper may be used for restricted orientation convex sets.

\section*{Acknowledgment}
The first author acknowledges Ho Chi Minh City University of Technology (HCMUT), VNU-HCM for  supporting this study.

\section*{Conflict of Interest Statement}
The authors have no conflicts of interest to declare.
All co-authors have seen and agree with the contents of the manuscript and there is no financial interest to report. We certify that this manuscript has not been submitted to, nor is under review at, another journal or other publishing venue.

\end{document}